\documentclass{article}
\usepackage{latexsym}
\usepackage{graphicx}
\usepackage{color}
\normalmarginpar
\title {\textbf{Shifted Chebyshev polynomials for Solving Three-Dimensional Volterra Integral Equations of the second kind}}
\author{ D.Sh.Mohamed*}
\date{\emph{Mathematics Department, Faculty of Science,\\ Zagazig University, Zagazig, Egypt\\ E-mail: doaashokry@zu.edu.eg}}
\begin{document}
\maketitle
\newcommand{\lb}{\hspace*{0.35in}}
\newcommand{\lc}{\hspace*{0.5in}}
\newcommand{\ld}{\hspace*{0.65in}}
\newcommand{\lf}{\hspace*{.80in}}
\newcommand{\lh}{\hspace*{.95in}}
\newcommand{\lj}{\hspace*{1.10in}}
\newcommand{\lk}{\hspace*{1.25in}}
\newcommand{\lm}{\hspace*{1.40in}}
\newcommand{\lo}{\hspace*{1.55in}}
\newcommand{\lp}{\hspace*{1.70in}}
\newcommand{\lr}{\hspace*{1.85in}}
\newcommand{\ls}{\hspace*{2.0in}}
\newcommand{\lt}{\hspace*{2.15in}}
\newcommand{\lv}{\hspace*{2.30in}}
\newcommand{\lx}{\hspace*{2.45in}}
\newcommand{\ly}{\hspace*{2.60in}}
\newcommand{\lz}{\hspace*{2.75in}}
\newcommand{\be}{\begin{equation}}
\newcommand{\ee}{\end{equation}}
\newcommand{\ba}{\begin{array}}
\newcommand{\ea}{\end{array}}
\newcommand{\goto}{\rightarrow}
\newcommand{\longto}{\longrightarrow}
\newcommand{\ov}{\overline}
\begin{abstract}
In this paper, an efficient method is presented for solving three
dimensional Volterra integral equations of the second kind with
continuous kernel. Shifted Chebyshev polynomial is applied to
approximate a solution for these integral equations. This method
transforms the integral equation to algebraic equations with
unknown Chebyshev coefficients. Numerical results are calculated
and the estimated error in each example is computed using Maple
17.

\end{abstract}
\textbf{keywords:} Three - dimensional Integral Equations, Shifted Chebyshev polynomials, collocation points. \\\\
\textbf{Mathematics Subject Classification: }  65R20, 41A50,
41A55, 65M70.
\section*{1. Introduction}
$~~~~~$Two dimensional integral equations provide an important
tool for modeling a numerous problems in engineering and science
[9,14]. There are many different numerical methods for solving one
and dimensional integral equations, such as
[1,3,4,5,7,9,10,13,15,16,17,18]. Some of these methods
can be used for solving three dimensional integral equations.\\\\
In [2,12] three dimensional differential transform method applied
for solving nonlinear three-dimensional Volterra integral
equations. The author in [11] applied Degenerate Kernel Method for
Three Dimension Nonlinear Integral Equations of the Second Kind.
The authors in [6] applied A computational method for nonlinear
mixed
Volterra-Fredholm integral equations.\\\\
In [17,18] the authors have been studied the numerical solutions
of two-dimensional integral equations by using Chebyshev
polynomials. The aim of the present paper is to solve three
dimensional Volterra integral equations of the second kind by
using Shifted Chebyshev polynomials. So, we consider
three-dimensional Volterra integral equations of the form:
$$u(x,y,z)=f(x,y,z)+\int_0^z \int_0^y \int_0^x k(x,y,z,r,s,t)u(r,s,t) dr ds dt,
\eqno(1.1)$$\\
where $(x,y,z) \in D=[0,X]\times [0,Y]\times [0,Z], u(x,y,z)$ is
the unknown function, $f(x,y,z)$ and $k(x,y,z,r,s,t)$ are given
functions defined, respectively on $D$ and\\
 $$G=\{ (x,y,z,r,s,t):
 0 \leq r\leq x \leq X,~0 \leq s\leq y \leq Y,~0 \leq s\leq z \leq Z\}.$$
\section*{2. Solution of three dimensional integral equations}
$~~~~~$ In this section Shifted Chebyshev polynomials  of the
first kind applied for Solving three dimensional Volterra integral
equations.\\\\
 This method is based on approximating the unknown function $u(x,y,z)$ as :\\\\
$$u(x,y,z)=\sum_{i=0}^\infty ~\sum_{j=0}^\infty ~\sum_{k=0}^\infty~a_{i,j,k}~ T^{*}_{i}(x)~ T^{*}_{j}(y)~T^{*}_{k}(z),~~~~~~
(x,y,z)\in D,~~~~~~~~~~~\eqno(2.1)$$\\
where $a_{i,j,k} , $  are constants to be determined, \\\\
 $T^{*}_{n} (x)$
is shifted Chebyshev polynomial of the first kind which is defined as [8]:\\
$$T^{*}_{n} (x) =T_{n} (2x-1)~,~~~T_{n} (x) =\cos n \theta~,~~ x=\cos \theta
,$$\\
and the following recurrence formulas :\\\\
$$T^{*}_{n} (x)=2( 2x-1) T^{*}_{n-1} (x)-T^{*}_{n-2} (x) ,~~~~~n=2,3,\ldots ,$$
with the initial conditions :\\ $$T^{*}_{0} (x)=1~,~~~~~~~~T^{*}_{1} (x)=2x-1 .$$\\\\
If the infinite series in (2.1) is truncated, then (2.1) can be written as :\\\\
$$ u(x,y,z)\approx u_N (x,y,z) =\sum_{i=0}^N ~\sum_{j=0}^N ~\sum_{k=0}^N~a_{i,j,k}~ T^{*}_{i}(x)~ T^{*}_{j}(y)~T^{*}_{k}(z),\eqno(2.2)$$\\
where $N$ is any natural number.\\\\
Substituting  from equation (2.2) into equation (1.1) we obtain :\\\\
$$\sum_{i=0}^N ~\sum_{j=0}^N ~\sum_{k=0}^N~a_{i,j,k}~~[ T^{*}_{i}(x)~ T^{*}_{j}(y)~T^{*}_{k}(z)
 - \int_0^z \int_0^y \int_0^x k(x,y,z,r,s,t)~~T^{*}_{i}(r)~ T^{*}_{j}(s)~T^{*}_{k}(t) dr ds dt]= f(x,y,z)~.\eqno(2.3)$$ \\
\\Hence the residual equation is defined as
:$$ R_N (x_{l},y_{m},z_{n})=\sum_{i=0}^N ~\sum_{j=0}^N
~\sum_{k=0}^N~a_{i,j,k}~~[~T^{*}_{i}(x)~
T^{*}_{j}(y)~T^{*}_{k}(z)-~~~~~~~~~~~~~~~~~~~~~~~~~~~~~~~~~~~~~~~~~~~~~~~~~$$
$$
 ~~~~~~~~~~~~~~~~~~~~~~~ \int_0^{z_{n}} \int_0^{y_{m}} \int_0^{x_{l}}~ k(x_l,y_m,z_n,r,s,t)~~T^{*}_{i}(r)~
  T^{*}_{j}(s)~T^{*}_{k}(t) dr ds dt~]- f(x_l,y_m,z_n)=0~,
  \eqno(2.4)$$\\\\
 \\ for Gauss - Chebyshev - lobatto collocation points [8]\\\\
  $$x_l
  =\frac{1}{2}[1+\cos(\frac{l\pi}{N})]~,~~~~~~~~~l=0,1,\ldots,N,$$\\
 $$y_m
  =\frac{1}{2}[1+\cos(\frac{m\pi}{N})]~,~~~~~~~~~m=0,1,\ldots,N,\eqno(2.5)$$\\
   $$z_n
  =\frac{1}{2}[1+\cos(\frac{n\pi}{N})]~,~~~~~~~~~n=0,1,\ldots,N.$$
\\
Equation (2.4)can be written as :\\\\
$$ \sum_{i=0}^N ~\sum_{j=0}^N
~\sum_{k=0}^N~a_{i,j,k}~~[~T^{*}_{i}(x)~
T^{*}_{j}(y)~T^{*}_{k}(z)-~~~~~~~~~~~~~~~~~~~~~~~~~~~~~~~~~~~~~~~~~~~~~~~~~~~~~~~~~~~~~~~~~$$$$
~~~~~~~~~~~~~~~~~~~~~~ \int_0^{z_{n}} \int_0^{y_{m}}
\int_0^{x_{l}}~ k(x_l,y_m,z_n,r,s,t)~~T^{*}_{i}(r)~
  T^{*}_{j}(s)~T^{*}_{k}(t) dr ds dt~]= f(x_l,y_m,z_n)~.
  \eqno(2.6)$$\\\\
Clearly, the obtained system of linear algebraic equations
contains $(N +1)^3$ equations in the same number as unknowns.
Solving this system we obtain the value of the constants
$a_{i,j,k}$ such that $i, j,k = 0 ,\ldots, N.$
\section*{3. Numerical Examples }
$~~~~$In this section some numerical examples of linear and
nonlinear three dimensional volterra integral equation are
 presented to illustrate the method.\\\\
 \textbf{Example 3.1} [12]:\\\\
Consider the following three dimensional Volterra integral equation :\\\\
$$u(x,y,z)=f(x,y,z)-\int_0^z \int_0^y \int_0^x u(r,s,t) dr ds dt,
\eqno(3.1)$$\\where $(x,y,z) \in [0,1]\times[0,1]\times[0,1]$
and\\
$$f(x,y,z)=x+y+z+\frac{x^2 y z+x y^2 z+xy z^2}{2}$$
with the exact solution $u(x,y,z)=x+y+z.$\\\\
Applying shifted Chebyshev polynomial of the first kind for
equation (3.1) when $N =1,$ we obtain\\\\\\\\\\\\
$$a_{0,0,0}+a_{0,0,1}(2z_n-1)+a_{0, 1, 0}(2y_m-1)+a_{0, 1, 1}(2z_n-1)(2y_m-1)+
a_{1, 0, 0}(2x_l-1)+a_{1, 0, 1}(2x_l-1)~(2z_n-1)+$$
$$a_{1, 1,
0}(2x_l-1)(2y_m-1)+a_{1, 1, 1} (2x_l-1)(2y_m-1)(2z_n-1)+
\frac{1}{2}(\frac{1}{2}(4a_{1, 0, 1}-4a_{1, 1, 1})x_l^2y_m-2a_{0,
1, 1} x_l y_m+$$
$$2a_{0, 0, 1}x_l y_m-2a_{1, 0,
1}x_l y_m +\frac{1}{2}(4a_{0, 1, 1}x_l+4a_{1, 1, 1} x_l^2-4a_{1,
1, 1}x_l)y_m^2+2a_{1, 1, 1}x_l y_m)z_n^2+$$
$$\frac{1}{2} (\frac{1}{2} (4 a_{1,1,0}-4 a_{1,1,1})
 x_l^2-2 a_{1,1,0} x_l-2 a_{0,1,1} x_l+2 a_{0,1,0} x_l+2 a_{1,1,1} x_l) y_m^2 z_n-a_{1,1,1} x_l y_m z_n+$$
$$\frac{1}{2}(2a_{1, 0, 0}-2a_{1, 0, 1}-2a_{1, 1, 0}+2a_{1, 1, 1})x_l^2y_m z_n-a_{0,0,1} x_l y_m z_n+a_{1,0,1} x_l y_m z_n-a_{0,1,0} x_l y_m z_n+$$
$$a_{0,0,0} x_l
y_m
 z_n-a_{1,0,0} x_l y_m z_n+a_{1,1,0} x_l y_m z_n+a_{0,1,1} x_l y_m z_n=x_l+y_m+z_n+\frac{1}{2}x_l^2y_m z_n+$$
$$~~~~~~~~~~~~~~~~\frac{1}{2}x_l y_m^2z_n+\frac{1}{2}x_l y_m z_n^2.\eqno(3.2)$$\\\\
Substituting  from the collocation points (2.5) when $N=1$ into
(3.2) we obtain a system of linear algebraic equations contains
$9$ equations in the same number as unknowns. Solving this system
we obtain the value of the constants as:
$$a_{0, 0, 0} = 3/2,~a_{0, 0, 1} = 1/2,~ a_{0, 1, 0} = 1/2,~
 a_{0, 1, 1} = 0,~ a_{1, 0, 0} = 1/2,~ a_{1, 0, 1} = 0,~ a_{1, 1, 0} = 0,~ a_{1, 1, 1} =
 0.$$\\
Substituting  from these constants into (2.2) we obtain the
approximate solution of equation (3.1) $u_1(x,y,z)=x+y+z$ which is
the same as the
exact solution.\\\\
 \textbf{Example 3.2} [12]:\\\\
Consider the following three dimensional Volterra integral equation :\\\\
$$u(x,y,z)=f(x,y,z)-24 x^2 y \int_0^z \int_0^y \int_0^x u(r,s,t) dr ds dt,
\eqno(3.3)$$\\where $(x,y,z) \in [0,1]\times[0,1]\times[0,1]$
and\\
$$f(x,y,z)=x^2 y+y z^2+x y z+24 x^2 y(\frac{1}{6} x y^2 z^3+\frac{1}{8} x^2y^2 z^2+\frac{1}{6} x^3 y^2 z)$$
with the exact solution $u(x,y,z)=x^2 y +y z^2+ x y z.$\\\\
Applying shifted Chebyshev polynomial of the first kind for
equation (3.3) when $N =2,$  and by using the collocations points
(2.5)we obtain a system of linear algebraic equations contains
$27$ equations in the same number as unknowns. Solving this system
we obtain the value of the constants as:\\\\
$$a_{0, 0, 0} = 1/2,~~a_{0, 0, 1} = 3/8,~~ a_{0, 0, 2} = 1/16,~~ a_{0, 1, 0} = 1/2,~~ a_{0, 1, 1} =
 3/8,~~ a_{0, 1, 2} = 1/16,~~ a_{0, 2, 0} = 0,$$
  $$ a_{0, 2, 1} =
 0,~~a_{0, 2, 2} = 0,~~a_{1, 0, 0} = 3/8,~~ a_{1, 0, 1} = 1/8,~~ a_{1, 0, 2} = 0,~~ a_{1, 1, 0} =
  3/8,~~
   a_{1, 1, 1} = 1/8,~~ a_{1, 1, 2} = 0,$$
   $$ a_{1, 2, 0} = 0,~~ a_{1, 2, 1} = 0,~~ a_{1, 2, 2} = 0,~~a_{2, 0, 0} = 1/16,~~ a_{2, 0, 1} = 0,~~ a_{2, 0, 2} = 0,~~ a_{2, 1, 0} =
   1/16,~~
   $$
   $$a_{2, 1, 1} = 0,~~ a_{2, 1, 2} = 0,~~ a_{2, 2, 0} = 0,~~ a_{2, 2, 1} = 0,~~ a_{2, 2, 2} =
   0.$$\\
Substituting  from these constants into (2.2) we obtain the
approximate solution of equation (3.3)\\ $u_2(x,y,z)=x^2 y +y z^2+
x y z$ which is the same as the
exact solution.\\\\
 \textbf{Example 3.3} [2]:\\\\
Consider the following nonlinear three dimensional Volterra integral equation :\\\\
$$u(x,y,z)=xyz-\frac{(xyz)^3}{27} \int_0^z \int_0^y \int_0^x u^2(r,s,t) dr ds dt,
\eqno(3.4)$$\\where $(x,y,z) \in [0,1]\times[0,1]\times[0,1]$
with the exact solution $u(x,y,z)=xyz.$\\\\
Applying shifted Chebyshev polynomial of the first kind for
equation (3.4) when $N =1,$  and by using the collocations points
(2.5)we obtain a system of nonlinear algebraic equations. Solving
this system
we obtain the value of the constants as:\\
$$a_{0, 0, 0} = 1/8,~ a_{0, 0, 1} = 1/8,~ a_{0, 1, 0} = 1/8,~ a_{0, 1, 1} = 1/8,~ a_{1, 0, 0} =
1/8,~
 a_{1, 0, 1} = 1/8,~ a_{1, 1, 0} = 1/8,~ a_{1, 1, 1} = 1/8.$$\\
Substituting  from these constants into (2.2) we obtain the
approximate solution of equation (3.4)\\ $u_1(x,y,z)=xyz$ which is
the same as the
exact solution.\\\\
\textbf{Example 3.4} [2]:\\\\
Consider the following linear three dimensional Volterra integral equation :\\\\
$$u(x,y,z)=x \cos z-\frac{(x^3y^3)}{9}\sin z+ \int_0^z \int_0^y \int_0^x r s^2~u(r,s,t) dr ds dt,
\eqno(3.5)$$\\where $(x,y,z) \in [0,1]\times[0,1]\times[0,1]$
with the exact solution $u(x,y,z)=x \cos z.$\\\\
Applying shifted Chebyshev polynomial of the first kind for
equation (3.5) when $N =2.$ Table 1 shows the absolute errors at
some particular points.\\\\
 \textbf{Table 1.Numerical results of
example 3.4 }
\begin{center}
\begin{tabular}{c c c c c}
  \hline\\
 $ x$ & $y$ & $z$ & Exact & Absolute error \\\\
  \hline\\\\
  0.1 & 0.1 & 0.1~~~ & 0.09950041653~~~~ & $2.2992069\times 10^{-4}$ \\\\\\
  0.01 & 0.1 & 0.1~~~ & 0.009950041653~~~~ & $2.2994289 \times 10^{-5}$
  \\\\\\
  0.01 & 0.01 & 0.1~~~ & 0.009950041653~~~~ & $2.2987394\times10^{-5}$ \\\\\\
  0.01 & 0.01 & 0.01~~~ & 0.009999500004~~~~ & $2.927109\times10^{-6}$ \\\\\\
  0.001 & 0.01 & 0.01~~~ & 0.0009999500004~~~~ & $2.927257\times10^{-7}$ \\\\\\
  0.001 & 0.001 & 0.01~~~ & 0.0009999500004~~~~ & $2.927145\times 10^{-7}$ \\\\\\
  0.001 & 0.001 & 0.001~~~ & 0.0009999995000~~~~ & $2.99266\times 10^{-8}$ \\\\
  \hline
\end{tabular}
\end{center}
\textbf{Example 3.5} [12]:\\\\
Consider the following three dimensional Volterra integral equation :\\\\
$$u(x,y,z)=f(x,y,z)+\int_0^z \int_0^y \int_0^x u(r,s,t) dr ds dt,
\eqno(3.6)$$
\\where $(x,y,z) \in [0,1]\times[0,1]\times[0,1]$
and\\
$$f(x,y,z)=e^{x+y}+e^{x+z}+e^{y+z}-e^x-e^y-e^z+1$$
with the exact solution $u(x,y,z)=e^{x+y+z}.$\\\\
Applying shifted Chebyshev polynomial of the first kind for
equation (3.6) when $N =2.$ Table 2 shows the absolute errors at
some particular points.\\\\\\
\textbf{Table 2.Numerical results of example 3.5 }\\\\
\begin{center}
\begin{tabular}{c c c c c}
  \hline\\
 $ x$ & $y$ & $z$ & Exact & Absolute error \\\\
  \hline\\\\
  0.1 & 0.1 & 0.1~~~ & 1.349858808~~~~ & $0.033089467$ \\\\\\
  0.01 & 0.1 & 0.1~~~ & 1.233678060~~~~ & 0.021664584
  \\\\\\
  0.01 & 0.01 & 0.1~~~ & 1.127496852~~~~ & 0.011933284 \\\\\\
  0.01 & 0.01 & 0.01~~~ & 1.030454534~~~~ & 0.003668302 \\\\\\
  0.001 & 0.01 & 0.01~~~ & 1.021222052~~~~ & 0.002550285 \\\\\\
  0.001 & 0.001 & 0.01~~~ & 1.012072289~~~~ & 0.001450893\\\\\\
  0.001 & 0.001 & 0.001~~~ & 1.003004505~~~~ & 0.000369862 \\\\
  \hline
\end{tabular}
\end{center}
 \textbf{Conclusion}\\
Analytical solution of the two and three dimensional integral
equations are usually difficult. In many cases, it is required to
approximate solutions. In this work, the three dimensional linear
and nonlinear integral equations of the second kind are solved by
using shifted Chebyshev polynomials through collocation points.

\end{document}